\documentclass[11pt,a4paper]{article}

\usepackage{amsmath,amsfonts,amssymb,amsthm}

\usepackage{bbm,mathrsfs, bm, mathtools}

\title{Threshold phenomena for random cones}

\author{Daniel Hug and Rolf Schneider}
\date{}
\newcommand{\Sd}{{\mathbb S}^{d-1}}
\newcommand{\R}{{\mathbb R}}

\newcommand{\C}{{\mathcal C}}

\newcommand{\bP}{{\mathbb P}}

\newcommand{\N}{{\mathbb N}}

\newcommand{\Ha}{\mathcal{H}}

\newcommand{\B}{\mathcal{B}}
\newcommand{\D}{{\rm d}}
\newcommand{\F}{{\mathcal F}}

\newcommand{\bE}{{\mathbb E}\,}
\newcommand{\bV}{{\mathbb V}\,}

\newcommand{\kuppa}{\kappa}

  \renewcommand{\exp}{{\rm exp}\,}

\newtheorem{theorem}{Theorem}
\newtheorem{lemma}{Lemma}

\begin{document}
\maketitle

\begin{abstract}
We consider an even probability distribution on the $d$-dimensional Euclidean space with the property that it assigns measure zero to any hyperplane through the origin. Given $N$ independent random vectors with this distribution, under the condition that they do not positively span the whole space, the positive hull of these vectors is a random polyhedral cone (and its intersection with the unit sphere is a random spherical polytope). It was first studied by Cover and Efron. We consider the expected face numbers of these random cones and describe a threshold phenomenon when the dimension $d$ and the number $N$ of random vectors tend to infinity. In a similar way, we treat the solid angle, and more generally the Grassmann angles. We further consider the expected numbers of $k$-faces and of Grassmann angles of index $d-k$ when also $k$ tends to infinity.\\[1mm]
\end{abstract}

\section{Introduction}\label{sec1}

The following is a literal quotation from \cite {DT09a}: ``Recent work has exposed a phenomenon of abrupt {\em phase transitions} in high-dimensional geometry. The phase transitions amount to a rapid shift in the likelihood of a property's occurrence when a dimension parameter crosses a critical level (a {\em threshold}).'' Two early observations in high-dimensional random geometry of such distinctly different behavior below and above a threshold were published in 1992. Dyer, F\"uredi and McDiarmid \cite{DFD92} considered the convex hull of $N=N(d)$ points chosen independently at random (with equal chances) from the vertices of the unit cube in $\R^d$. Let $V_{d,N}$ denote the volume of this random polytope. Then, for every $\varepsilon>0$,
$$ \lim_{d\to\infty} \bE V_{d,N} = \left\{\begin{array}{ll} 1, & \mbox{if }N\ge(2e^{-1/2}+\varepsilon)^d,\\[1mm]
0, & \mbox{if }N\le(2e^{-1/2}-\varepsilon)^d.\end{array}\right.$$
Here $\bE$ denotes mathematical expectation. The paper \cite{DFD92} has a similar result for the convex hull of i.i.d.~uniform random points from the interior of the unit cube. We have quoted this example as an illustration of what we have in mind: for instance, a $d$-dimensional random polytope with its number $N$ of vertices depending on $d$, where a small change of this dependence causes an abrupt change of some property as $d\to\infty$. In the work of Vershik and Sporyshev \cite{VS92}, a $d$-dimensional random polytope is obtained as a uniform random orthogonal projection of a fixed regular simplex with $N$ vertices in a higher-dimensional space, and threshold phenomena are exhibited for the expected numbers of $k$-faces, under the assumption of a linearly coordinated growth of the parameters $d,N,k$. Similar models, also with the regular simplex replaced by the regular cross-polytope, and random projections extended to more general random linear mappings, have found important applications in the work of Donoho and collaborators. We refer to the article of Donoho and Tanner \cite{DT09b}, where also earlier work of these authors is cited and explained. The paper \cite{DT10} of the same authors treats random projections of the cube and the positive orthant in a similar way. Generally in stochastic geometry, threshold phenomena have been investigated for face numbers, neighborliness properties, volumes, intrinsic volumes, more general measures, and for several different models of random polytopes. Different phase transitions were exhibited. We mention that \cite{GG09} has extended the model of \cite{DFD92} by introducing more general distributions for the random points. The paper \cite{Piv07} considers convex hulls of i.i.d.~random points with either Gaussian distribution or uniform distribution on the unit sphere. In \cite{BCGTT19}, \cite{BKT19}, the points have a beta or beta-prime distribution. The paper \cite{BO19} studies facet numbers of convex hulls of random points on the unit sphere in different regimes. The papers \cite{Piv07}, \cite{HHRT15} and \cite{BCGTT19} deal also with polytopes generated by intersections of random closed halfspaces. An important role is played by phase transitions in convex programs with random data. We quote from the work of Amelunxen, Lotz, McCoy and Tropp \cite{ALMT14}, which discovers and describes several of these phenomena: ``This paper provides the first rigorous analysis that explains why phase transitions are ubiquitous in random convex optimization problems. \dots The applied results depend on foundational research in conic geometry.''

In this paper, we consider a model of random polyhedral convex cones (or, equivalently, of random spherical polytopes) that was introduced by Cover and Efron \cite{CE67} (and more closely investigated in \cite{HS16}). Let $\phi$ be a probability measure on the Euclidean space $\R^d$ which is even (invariant under reflection in the origin $o$) and assigns measure zero to each hyperplane through the origin. For $n\in\N$, the {\em $(\phi,n)$-Cover--Efron cone} $C_n$ is defined as the positive hull of $n$ independent random vectors $X_1,\dots,X_n$ with distribution $\phi$, under the condition that this positive hull is different from $\R^d$. The intersection $C_n\cap \Sd$ with the unit sphere $\Sd$ is a spherical random polytope, contained in some closed hemisphere. In the following it will be convenient to work with polyhedral cones instead of spherical polytopes.

For $k\in\{1,\dots,d-1\}$, let $f_k(C_n)$ denote the number of $k$-dimensional faces of the cone $C_n$ (equivalently, the number of $(k-1)$-dimensional faces of the spherical polytope $C_n\cap \Sd$). We are interested in the asymptotic behavior of the expectation $\bE f_k(C_n)$, as $d$ tends to infinity and $n$ grows suitably with $d$.

\vspace{2mm}

\noindent{\bf Convention.} In the following, $C_N$ is a $(\phi,N)$-Cover--Efron cone in $\R^d$, and $N=N(d)\ge d$ is an integer depending on the dimension $d$, but we will omit the dimension $d$ in the notation.

\vspace{2mm}

Any $k$-face of the cone $C_N$ is a.s. the positive hull of $k$ vectors from $X_1,\dots,X_N$, and there are $\binom{N}{k}$ possible choices. Therefore, we consider the quotient $\bE f_k(C_N)/\binom{N}{k}$. For this, we can state the following behavior above and below a threshold.

\begin{theorem}\label{T1.1}
Suppose that
$$ \frac{d}{N} \to\delta \quad\mbox{as }d\to\infty,$$
with a number $\delta\in[0,1]$. Let $k\in\N$ be fixed. Then the Cover--Efron cone $C_N$ satisfies
$$ \lim_{d\to\infty} \frac{\bE f_k(C_N)}{\binom{N}{k}} = \left\{\begin{array}{ll} 1 & \mbox{if } 1/2<\delta\le1,\\[1mm]
(2\delta)^k & \mbox{if } 0\le \delta< 1/2.\end{array}\right.$$
\end{theorem}

We have assumed here that the number $k$ is constant. But if in the case where $0\le \delta< 1/2$ we let also $k$ depend on $d$, in such a way that
$$ k\to\infty \quad \mbox{and}\quad \frac{k}{d}\to 0\quad\mbox{as }d\to\infty,$$
(thus, $k$ grows with $d$, but sublinearly) then
$$ \lim_{d\to\infty} \frac{\bE f_k(C_N)}{\binom{N}{k}} = 0.$$
This is easily read off from the proof of the second part of Theorem \ref{T1.1} given in Section \ref{sec4}.

The case $1/2<\delta\le 1$ of Theorem \ref{T1.1} follows from the stronger Theorem \ref{T1.2a} below, but we have formulated Theorem \ref{T1.1} in this way to illuminate the different behavior if $\delta$ is below $1/2$ or above $1/2$. Partial information on the case $\delta= 1/2$ is contained in Theorem \ref{T1.2a}. In this theorem, also the number $k$ is allowed to depend on $d$.

\begin{theorem}\label{T1.2a}
Let $k=k(d)\in \{1,\dots,d-1\}$. Suppose that
$$\frac{k}{\sqrt{N}}\to 0 \mbox{ as } d\to\infty\qquad \text{and}\qquad \sqrt{N}\left(\frac{d}{N}-\frac{1}{2}\right)\ge  a $$
with some real constant $a$. Then
$$ \lim_{d\to\infty} \frac{\bE f_k(C_N)}{\binom{N}{k}} =1. $$
The same conclusion is obtained, if $k=k(d)$ is bounded (as $d\to\infty$) and
$$
\sqrt{N}\left(\frac{d}{N}-\frac{1}{2}\right)\ge - a\,\left(\log\log N\right)^{\alpha}
$$
with  constants $a>0$ and $\alpha\ge 0$.
\end{theorem}

The assumptions of this theorem are, in particular, satisfied if $k$ is constant and
$$
N-2d\le 2a\sqrt{N}\left(\log\log N\right)^{\alpha},
$$
(which holds, for instance, if $N-2d$ is bounded from above); the latter holds if $\frac{d}{N}\to\delta>\frac{1}{2}$.

We do not have complete information in the case where $\frac{d}{N}\to\frac{1}{2}$, but there is a precise result if $\frac{d}{N}=\frac{1}{2}$.

\begin{theorem}\label{T1.2b}
If $N=2d$ and $k\in \N$ is fixed, then
$$ \lim_{d\to\infty} \sqrt{d}\left(1-\frac{\bE f_k(C_N)}{\binom{N}{k}}\right)=\frac{k}{\sqrt{\pi}}.$$
\end{theorem}

On the other hand, Theorem \ref{T1.5} below implies that if
$$ \frac{d}{N} \to\delta\in \left(0,\tfrac{1}{2}\right] \qquad \text{and}\qquad  \frac{k}{d}\to \rho \in (0,1) \quad\mbox{as }d\to\infty,$$
then
$$
\lim_{d\to\infty} \frac{\bE f_k(C_N)}{\binom{N}{k}} = 0.
$$

As another functional of a closed convex cone $C\subset\R^d$, we consider the solid angle $v_d(C)$. This is the normalized spherical Lebesgue measure of $C\cap\Sd$. (We avoid the notation $V_d(C)$ used in \cite{HS16}, since $V_d$ is often used for the volume of a convex body. The reader is warned that what we denote here by $v_d(C)$ was denoted by $v_{d-1}(C\cap\Sd)$ in \cite[Sect.~6.5]{SW08}.)

More generally, we consider the Grassmann angles. For a closed convex cone $C\subset\R^d$ which is not a subspace, the $j$th {\em Grassmann angle} of $C$, for $j\in\{1,\dots,d-1\}$, is defined by
$$ U_j(C):= \frac{1}{2}\bP(C\cap{\mathcal L}\not=\{o\}\},$$
where ${\mathcal L}$ is a random $(d-j)$-dimensional subspace with distribution $\nu_{d-j}$. The latter is the unique Haar probability measure on $G(d,d-j)$, the Grassmannian of $(d-j)$-dimensional linear subspaces of $\R^d$. Thus,
$$ U_j(C)= \frac{1}{2} \int_{G(d,d-j)} {\mathbbm 1}\{C\cap L\not=\{o\}\}\,\nu_{d-j}(\D L).$$
Grassmann angles were introduced by Gr\"unbaum \cite{Gru68}, in a slightly different, though equivalent way. Gr\"unbaum's Grassmann angles are given by $\gamma^{d-j,d}=1- 2U_j$. We note that $v_d=U_{d-1}$, and that $U_j(C)\le 1/2$, with equality if $C$ is a halfspace.

\begin{theorem}\label{T1.3}
Suppose that
$$ \frac{d}{N} \to\delta \quad\mbox{as }d\to\infty,$$
with a number $\delta\in[0,1]$. Let $k\in\N$ be fixed. Then the $(\phi,N)$-Cover--Efron cone $C_N$ satisfies
$$ \lim_{d\to\infty} \bE2 U_{d-k}(C_N) = \left\{\begin{array}{ll} 0 &  \mbox{if } 1/2<\delta \le 1,\\[1mm]
1-\left(\frac{\delta}{1-\delta}\right)^k & \mbox{if }  0\le \delta< 1/2.\end{array}\right.$$
\end{theorem}

We have assumed here that the number $k$ is constant. But if in the case where $0\le \delta< 1/2$ we let also $k$ depend on $d$, in such a way that
$$ k\to\infty \quad \mbox{and}\quad \frac{k}{d}\to 0\quad\mbox{as }d\to\infty,$$
(thus, $k$ grows with $d$, but sublinearly) then
$$ \lim_{d\to\infty}  \bE 2 U_{d-k}(C_N) = 1.$$
This can be seen from the second part of the proof of Theorem \ref{T1.3} given in Section \ref{sec5}.

Again, the first part of this theorem has a stronger version, given by the following theorem.

\begin{theorem}\label{T1.4a}
Let $k=k(d)\in \{1,\dots,d-1\}$. Suppose that
$$\frac{k}{\sqrt{N}}\to 0 \mbox{ as } d\to\infty\qquad \text{and}\qquad \sqrt{N}\left(\frac{d}{N}-\frac{1}{2}\right)\ge  a $$
with some real constant $a$. Then
$$ \lim_{d\to\infty} \bE U_{d-k}(C_N)=0. $$
The same conclusion is obtained, if $k=k(d)$ is bounded (as $d\to\infty$) and
$$
\sqrt{N}\left(\frac{d}{N}-\frac{1}{2}\right)\ge - a\,\left(\log\log N\right)^{\alpha}
$$
with  constants $a>0$ and $\alpha\ge 0$.
\end{theorem}

And similarly as above, there is a more precise asymptotic relation if $N=2d$.

\begin{theorem}\label{T1.4b}
If $N=2d$ and $k\in\N$ is fixed, then
$$ \lim_{d\to\infty} \sqrt{d}\cdot\bE U_{d-k}(C_N) =\frac{k}{\sqrt{\pi}}.$$
\end{theorem}

On the other hand, Theorem \ref{T1.7} below implies that if
$$ \frac{d}{N} \to\delta\in \left(0,\tfrac{1}{2}\right] \qquad \text{and}\qquad  \frac{k}{d}\to \rho \in (0,1) \quad\mbox{as }d\to\infty,$$
then
$$
\lim_{d\to\infty}  \bE 2U_{d-k}(C_N) = 1.
$$

Clearly, in Theorem \ref{T1.1} (and similarly in Theorem \ref{T1.3}), the change when passing a threshold is not so abrupt as in the examples from \cite{DT10}, where both parameters, $N$ and $k$, grow linearly with the dimension:  below the threshold $\delta=1/2$, the limit in question increases (decreases) with the parameter to an extremal value, above the threshold, it remains constant. The situation changes if also the number $k$ increases sublinearly in the dimension. Then indeed we also have a sharp threshold as pointed out above.  Now we consider the case where $k$ increases proportional to the dimension; then a more subtle phase transition is observed. Under a linearly coordinated growth, for $k$-faces we find the same threshold as established by Donoho and Tanner \cite{DT10} in their investigation of random linear images of orthants. This may seem unexpected, since the random cones considered in \cite{DT10} and here have different distributions (see, however, the Appendix).

As in \cite{DT10}, we define
$$ \rho_W(\delta):= \max\{0,2-\delta^{-1}\} \quad\mbox{for }0<\delta<1$$
(the index $W$ stands for `weak' threshold).

\begin{theorem}\label{T1.5}
Let $0<\delta<1$ and $0\le \rho<1$ be given. Let $k(d)=k<d<N=N(d)$ be integers such that
$$ \frac{d}{N} \to\delta, \qquad \frac{k}{d}\to \rho \qquad\mbox{as }d\to\infty.$$
Then
$$ \lim_{d\to\infty} \frac{\bE f_k(C_N)}{\binom{N}{k}} =\left\{ \begin{array}{ll} 1 & \mbox{if }\rho<\rho_W(\delta),\\ 0 & \mbox{if }\rho >\rho_W(\delta).\end{array}\right.$$
\end{theorem}

We note that the first assumption of this theorem, $\rho<\rho_W(\delta)$, implies that for large $d$ we have $N+k<2d$.

Adapting an argument of Donoho and Tanner \cite{DT10} to the present situation, we can also replace the convergence of an expectation in the first part of Theorem \ref{T1.5} by the convergence of a probability, at the cost of a smaller threshold.

\begin{theorem}\label{T1.6}
Let $0<\delta,\rho<1$ be given, where $\delta>1/2$. Let $k(d)=k<d<N=N(d)$ be integers such that
$$ \frac{d}{N}\to\delta,\qquad \frac{k}{d}\to\rho\qquad\mbox{as } d\to\infty.$$
Then there exists a positive number $\rho_S(\delta)$ such that, for $\rho<\rho_S(\delta)$,
$$ \lim_{d\to\infty} \left[\bE f_k(C_N)-\binom{N}{k}\right]=0\quad\mbox{and}\quad\lim_{d\to\infty} \bP\left(f_k(C_N) =\binom{N}{k}\right)=1.$$
\end{theorem}

There is also a counterpart to Theorem \ref{T1.5} for Grassmann angles.

\begin{theorem}\label{T1.7}
Let $0<\delta<1$ and $0\le \rho<1$ be given. Let $k(d)=k<d<N=N(d)$ be integers such that
$$ \frac{d}{N} \to\delta, \qquad \frac{k}{d}\to \rho \qquad\mbox{as }d\to\infty.$$
Then
$$ \lim_{d\to\infty} \bE 2U_{d-k}(C_N)  =\left\{ \begin{array}{ll}  0 & \mbox{if }\rho <\frac{1}{2}\rho_W(\delta),\\[2mm]
1 & \mbox{if }\rho > \frac{1}{2}\rho_W(\delta).\end{array}\right.$$
\end{theorem}

After some preliminaries in the next section, we collect a number of auxiliary results about sums of binomial coefficients in Section \ref{sec3}. Then we prove the first three theorems in Section \ref{sec4}, Theorems \ref{T1.3} to \ref{T1.4b} in Section \ref{sec5}, Theorems \ref{T1.5} and \ref{T1.6} in Section \ref{sec6}, and Theorem \ref{T1.7} in Section \ref{sec7}.

\section{Preliminaries}\label{sec2}

First we recall two classical facts. For $n\in\N$, let $H_1,\dots,H_n\in G(d,d-1)$. Assume that these hyperplanes are in general position, that is, the intersection of any $m\le d$ of them is of dimension $d-m$. Then the number of $d$-dimensional cones in the tessellation of $\R^d$ induced by these hyperplanes is given by
$$ C(n,d) := 2\sum_{i=0}^{d-1}\binom{n-1}{i}.$$
From this result of Steiner (in dimension three) and Schl\"afli, Wendel has deduced the following. If $X_1,\dots,X_n$ are i.i.d.~random vectors in $\R^d$ with distribution $\phi$ (enjoying the properties mentioned above), then
$$ P_{d,n}:= \bP\left({\rm pos}\{X_1,\dots,X_n\}\not=\R^d\right) =\frac{C(n,d)}{2^n},$$
where $\bP$ stands for probability and ${\rm pos}$ denotes the positive hull. For references and proofs, we refer to \cite[Sect.~8.2.1]{SW08}. Now we can write down the distribution of the $(\phi,n)$-Cover--Efron cone $C_n$, namely
$$
\bP\left(C_n\in B\right)= \frac{1}{P_{d,n}} \int_{(\Sd)^n} {\mathbbm 1}\{{\rm pos}\{x_1,\dots,x_n\}\in B\setminus\{\R^d\}\}\,\phi^n(\D(x_1,\dots,x_n))
$$
for $B\in\B(\C^d)$, where $\C^d$ denotes the space of closed convex cones in $\R^d$ (with the topology of closed convergence) and $\B(\C^d)$ is its Borel $\sigma$-algebra.

There is an equivalent representation of $C_n$. For this, we denote by $\phi^*$ the image measure of $\phi$ under the mapping $x\mapsto x^\perp$ from $\R^d\setminus\{o\}$ to $G(d,d-1)$.
Let $\Ha_1,\dots,\Ha_n$ be i.i.d.~random hyperplanes with distribution $\phi^*$. They are almost surely in general position. The {\em $(\phi^*,n)$-Schl\"afli cone} $S_n$ is obtained by picking at random (with equal chances) one of the $d$-dimensional cones from the tessellation induced by $\Ha_1,\dots,\Ha_n$. Its distribution is given by
$$
\bP\left(S_n\in B\right) = \int_{G(d,d-1)^n}\frac{1}{C(n,d)}\sum_{C\in\F_d(H_1,\dots,H_n)} {\mathbbm 1}\{C\in B\}\,\phi^{*n}(\D(H_1,\dots,H_n))
$$
for $B\in\B(\C^d)$, where $\F_d(H_1,\dots,H_n)$ is the set of $d$-cones in the tessellation induced by $H_1,\dots,H_n$.
We have (see \cite[Thm.~3.1]{HS16})
$$
C_n= S_n^\circ\quad \mbox{in distribution},
$$
where $S_n^\circ$ denotes the polar cone of $S_n$.

For the expectations appearing in our theorems, explicit representations are available. The proofs  of Theorems \ref{T1.1}, \ref{T1.2a}, \ref{T1.2b} and \ref{T1.5} are based on the formula
\begin{equation}\label{3.3}
\frac{\bE f_k(C_N)}{\binom{N}{k}}= 2^k\frac{C(N-k,d-k)}{C(N,d)} =\frac{P_{d-k,N-k}}{P_{d,N}}
\end{equation}
for $k\in\{0,\ldots,d-1\}$ (see \cite[(3.3)]{CE67} or \cite[(27)]{HS16}). For the proofs of Theorems \ref{T1.3}, \ref{T1.4a}, \ref{T1.4b}  and \ref{T1.7}, we use the explicit formula
\begin{equation}\label{3.3a}
\bE 2 U_{d-k}(C_N) = \frac{C(N,d)-C(N,d-k)}{C(N,d)} = 1-\frac{P_{d-k,N}}{P_{d,N}}
\end{equation}
for $k\in\{1,\ldots,d-1\}$ (see \cite[(29)]{HS16}). It is sometimes useful to write this in the form
\begin{equation}\label{4.r}
\bE 2U_{d-k}(C_N) = \frac{1+\binom{N-1}{d-1}^{-1}\sum_{i=d-k}^{d-2}\binom{N-1}{i}} {1+\binom{N-1}{d-1}^{-1}\sum_{i=0}^{d-2}\binom{N-1}{i}}
\end{equation}
(where an empty sum is zero, by definition).

\section{Auxiliary results on binomial coefficients}\label{sec3}

First we collect some information on the Wendel probabilities
$$ P_{d,n} = \frac{1}{2^{n-1}}\sum_{i=0}^{d-1}\binom{n-1}{i}.$$

Let $\xi_n$ be a random variable which has the binomial distribution with parameters $n$ and $1/2$, thus
$$ \bP(\xi_n=k) =\frac{1}{2^n}\binom{n}{k} \quad\mbox{and}\quad \bP(\xi_n\le k) =\frac{1}{2^n}\sum_{i=0}^k \binom{n}{i}.$$
Thus we can write
\begin{equation}\label{2.1a}
P_{d,n} = \bP(\xi_{n-1}\le d-1)
\end{equation}
and therefore, by (\ref{3.3}),
\begin{equation}\label{2.1b}
\frac{ \bE f_k(C_N)}{\binom{N}{k}} =\frac{\bP(\xi_{N-k-1}\le d-k-1)}{\bP(\xi_{N-1}\le d-1)}.
\end{equation}
Similarly, by (\ref{3.3a}) we have
\begin{equation}\label{2.1c}
\bE 2U_{d-k}(C_N) =1-\frac{\bP(\xi_{N-1}\le d-k-1)}{\bP(\xi_{N-1}\le d-1)}.
\end{equation}

The following two lemmas concern the Wendel probabilities and are therefore stated here, although they are not needed before the proof of Theorem \ref{T1.6}.

\begin{lemma}\label{L3.1}
For $k\in\{1,\dots,d-1\}$,
\begin{equation}\label{3.1a}
\left(\frac{\bE f_k(C_N)}{\binom{N}{k}}\right)^{-1}=  \frac{P_{d,N}}{P_{d-k,N-k}} = 1+ A
\end{equation}
with
$$ A=\frac{1}{2^{N-1}P_{d-k,N-k}}\, {\sum_{j=1}^k \binom{k}{j}\sum_{m=0}^{j-1}\binom{N-k-1}{d-k+m}}.$$
\end{lemma}

\begin{proof}
Writing $\binom{N-1}{i}= \binom{N-2}{i-1}+\binom{N-2}{i}$ for $i=1,\dots,d-1$, we obtain
$$ P_{d,N}=\frac{1}{2}P_{d-1,N-1}+\frac{1}{2}P_{d,N-1}.$$
This and induction can be used to prove that for $k\in\{1,\dots,d-1\}$,
$$P_{d,N} =\frac{1}{2^k}\sum_{j=0}^k\binom{k}{j} P_{d-k+j,N-k}.$$
For $j\in\{1,\dots,k\}$ we have
$$ P_{d-k+j,N-k} = P_{d-k,N-k} +\frac{1}{2^{N-k-1}}\sum_{m=0}^{j-1}\binom{N-k-1}{d-k+m}.$$
This gives
\begin{eqnarray*}
\frac{P_{d,N}}{P_{d-k,N-k}} &=& \frac{1}{2^k}\sum_{j=0}^k\binom{k}{j}\frac{P_{d-k+j,N-k}}{P_{d-k,N-k}}\\
&=& \frac{1}{2^k}\left[ \binom{k}{0} +\sum_{j=1}^k\binom{k}{j}\left(1+\frac{1}{2^{N-k-1}P_{d-k,N-k}} \sum_{m=0}^{j-1}\binom{N-k-1}{d-k+m}\right)\right]
\end{eqnarray*}
and thus the assertion.
\end{proof}

\begin{lemma}\label{L3.2}
$$ P_{d-k,N-k}\ge \frac{1}{2}-\frac{1}{2^{N-k}}\sum_{r=0}^{N-2d+k-1}\binom{N-k-1}{d-k+r}.$$
If $N-2d+k-1<0$, the sum is zero, by convention.
\end{lemma}

\begin{proof}
Let $M,p$ be integers. If $1\le p\le\frac{M}{2}$, then
\begin{eqnarray*}
2^M = \sum_{i=0}^M \binom{M}{i} &=& \sum_{i=0}^{p-1}\binom{M}{i} +\binom{M}{p}+\dots+\binom{M}{M-p} +\sum_{i=M-p+1}^M
\binom{M}{i}\\
&=& 2\sum_{i=0}^{p-1} \binom{M}{i} + \sum_{r=0}^{M-2p}\binom{M}{p+r},
\end{eqnarray*}
thus
$$ \sum_{i=0}^{p-1}\binom{M}{i} = 2^{M-1} -\frac{1}{2}\sum_{r=0}^{M-2p}\binom{M}{p+r}.$$

If $p> \frac{M}{2}$, we have
$$ 2^M = \sum_{i=0}^M \binom{M}{i}\le 2\sum_{i=0}^{p-1}\binom{M}{i}.$$
Hence, for arbitrary $p\ge 1$ we may write
$$
\sum_{i=0}^{p-1}\binom{M}{i} \ge 2^{M-1} -\frac{1}{2}\sum_{r=0}^{M-2p}\binom{M}{p+r},
$$
with the convention that the last sum is zero if $2p> M$. The choice $M= N-k-1$ and $p=d-k$ now gives the assertion.
\end{proof}

Below some information on binomial coefficients is required. First we note Stirling's formula
\begin{equation}\label{Stirling}
 n!= \sqrt{2\pi n}\,e^{-n}n^ne^{\theta/12n},\quad 0<\theta<1.
\end{equation}
It implies, in particular, that
\begin{equation}\label{2.1}
\binom{2n}{n}\sim \frac{2^{2n}}{\sqrt{\pi n}} \quad\mbox{as }n\to\infty
\end{equation}
(where $a_n\sim b_n$ as $n\to\infty$ means that $\lim_{i\to\infty}a_n/b_n=1$).

The following lemma gives upper and lower bounds for the expressions appearing in (\ref{4.r}). For the proof of the upper bound, we adjust and slightly refine the argument for Proposition 1(c) in \cite{Klar}, in the current framework. The improved lower bound in \eqref{4.10} will be crucial in
the following.

\begin{lemma}\label{L4.1}
Let $n\in\N$, $m\in\N_0$ and $2m\le n+1$.
\begin{enumerate}
\item[{\rm (a)}]
If $2m<n+1$, then
\begin{equation}\label{4.11}
\frac{1}{\binom{n}{m+1}}\sum_{j=0}^{m}\binom{n}{j}
\le \frac{m+1}{n-m}\cdot\frac{n-m+1}{n-2m+1}\left(1-\left(\frac{m}{n-m+1}\right)^{m+1}\right).
\end{equation}
If $2m=n+1$, then
\begin{equation}\label{4.11b}
\frac{1}{\binom{n}{m+1}}\sum_{j=0}^{m}\binom{n}{j}
\le \frac{(m+1)^2}{n-m}.
\end{equation}
\item[{\rm (b)}]
If $2\le \ell \le m$, then
\begin{equation}\label{4.10}
\frac{m-\ell+1}{n-2m+2\ell-1}\left(1-\left(\frac{m-\ell+1}{n-m+\ell}\right)^{\ell+1}\right)
\le \frac{1}{\binom{n}{m+1}}\sum_{j=0}^{m}\binom{n}{j}.
\end{equation}
Moreover,
\begin{equation}\label{4.10neu}
\frac{m+1}{n-m}\le\frac{m+1}{n-m}\cdot\frac{n+1}{n+1-m}\le  \frac{1}{\binom{n}{m+1}}\sum_{j=0}^{m}\binom{n}{j}.
\end{equation}
\end{enumerate}
\end{lemma}

\begin{proof} (a) The cases $n=1$, $m\in\{0,1\}$ and $n\ge 2$, $m=0$ are easy to check.
Now let $n\ge 2$, $m\ge 1$, and hence also $m\le n-1$. If $j\in\{0,\ldots,m-1\}$, then
$$
\frac{\binom{n}{j}}{\binom{n}{m}}=\frac{m}{n-m+1}\cdots\frac{j+1}{n-j}
\le \left(\frac{m}{n-m+1}\right)^{m-j},
$$
since
$$
\frac{j+1}{n-j}=-1+\frac{n+1}{n-j}\le -1+\frac{n+1}{n-(m-1)}=\frac{m}{n-m+1}\in (0,1].
$$
Therefore, if $2m<n+1$, then $0<q_0:=m/(n-m+1)<1$ and
$$
\frac{1}{\binom{n}{m+1}}\sum_{j=0}^m\binom{n}{j}\le \frac{\binom{n}{m}}{\binom{n}{m+1}}
\sum_{j=0}^m q_0^{m-j}=\frac{m+1}{n-m}\cdot\frac{1-q_0^{m+1}}{1-q_0},
$$
which implies \eqref{4.11}. If $2m=n+1$, then $m/(n-m+1)=1$, and \eqref{4.11b} follows similarly.

(b)  Note that $m\le ({n+1})/{2}\le n$, and
$m\le n-1$ if $n\ge 2$. Hence, if $ 2\le \ell\le m$, then $m\le n-1$, $n-2m+2\ell -1\ge 2$,
$(n+1)/(m+1)>1$ and $0<q_1:=(m-\ell +1)/(n-m+\ell)<1$. Then, for $j\in\{m-\ell,\ldots,m\}$  we obtain
\begin{align*}
\frac{\binom{n}{j}}{\binom{n}{m+1}}
&=\left(\frac{n+1}{m+1}-1\right)^{-1}\cdots \left(\frac{n+1}{j+1}-1\right)^{-1}
\ge \left(\frac{n+1}{j+1}-1\right)^{-(m+1-j)}\\
&\ge \left(\frac{m-\ell+1}{n-m+\ell}\right)^{m+1-j}= q_1^{m+1-j}
\end{align*}
and hence
$$
\frac{1}{\binom{n}{m+1}}\sum_{j=m-\ell}^{m}\binom{n}{j} \ge q_1\sum_{r=0}^\ell q_1^{r}
= \frac{q_1}{1-q_1}\left(1-q_1^{\ell+1}\right),
$$
which yields \eqref{4.10}.

If $m=n$, then \eqref{4.10neu} holds trivially, since $\binom{n}{n+1}=0$. It also holds for $m=0$. In the remaining cases, we have
$$
\frac{1}{\binom{n}{m+1}}\sum_{j=0}^{m}\binom{n}{j}
\ge \frac{\binom{n}{m}+\binom{n}{m-1}}{\binom{n}{m+1}}=
\frac{m+1}{n-m}\cdot\frac{n+1}{n+1-m}.
$$
This completes the proof (b).
\end{proof}

From \eqref{4.11} and \eqref{4.10neu} with $n=N-1$ and $m=d-2$, we deduce that
\begin{equation}\label{4.2}
\frac{d-1}{N-d+1}\le \frac{1}{\binom{N-1}{d-1}}\sum_{j=0}^{d-2}\binom{N-1}{j}\le
\frac{d-1}{N-d+1}\cdot\frac{N-d+2}{N-2d+4},
\end{equation}
if $N>2d-4$.

\begin{lemma}\label{L4.4}
If $d/N \to\delta$ as $d\to\infty$, with $0\le \delta <1/2$, then
$$
\lim_{d\to\infty}
\frac{1}{\binom{N-1}{d-1}}\sum_{j=0}^{d-2}\binom{N-1}{j}= \frac{\delta}{1-2\delta}.
$$
\end{lemma}

\begin{proof}
Assume that $d/N\to \delta$ as $d\to\infty$, with a number $0\le \delta < 1/2$. We write $N=\alpha d$, where $\alpha$ depends on $d$ and satisfies $\alpha\to\delta^{-1}$ as $d\to\infty$. If $\delta=0$, this means that $\alpha\to\infty$. We assume that $d$ is so large that $\alpha>2$. From (\ref{4.2}) we have
$$ \frac{1}{\binom{N-1}{d-1}}\sum_{j=0}^{d-2}\binom{N-1}{j}\le \frac{d-1}{(\alpha-1)d+1}\cdot\frac{(\alpha-1)d+2}{(\alpha-2)d+4}.$$
We conclude that
\begin{equation}\label{4.6}
\limsup_{d\to\infty} \frac{1}{\binom{N-1}{d-1}}\sum_{j=0}^{d-2}\binom{N-1}{j}\le \frac{\delta}{1-2\delta}.
\end{equation}

Lemma \ref{L4.1}(b) provides the lower bound
$$ \frac{1}{\binom{N-1}{d-1}}\sum_{j=0}^{d-2} \binom{N-1}{j} \ge \frac{d-\ell_1-1}{N-2d+2\ell_1+2}\left(1-\left(\frac{d-\ell_1-1}{N-d+\ell_1+1}\right)^{\ell_1+1}\right),$$
if $2\le \ell_1\le d-2$. From this, we obtain
$$ \liminf_{d\to\infty}\frac{1}{\binom{N-1}{d-1}}\sum_{j=0}^{d-2}\binom{N-1}{j}
\ge \frac{\delta}{1-2\delta}\left(1-\left(\frac{\delta}{1-\delta}\right)^{\ell_1+1}\right)$$
for each fixed $\ell_1\ge 2$. Letting $\ell_1\to\infty$, we find that
$$
\liminf_{d\to\infty}
\frac{1}{\binom{N-1}{d-1}}\sum_{j=0}^{d-2}\binom{N-1}{j}\ge \frac{\delta}{1-2\delta}.
$$
Together with (\ref{4.6}) this completes the proof.
\end{proof}

We state another simple lemma.

\begin{lemma}\label{L4.2}
Let $m\in\N$. Then
$$
\sum_{i=0}^m \binom{2m}{i} = 2^{2m-1}+\frac{1}{2}\binom{2m}{m},\qquad \sum_{i=0}^{m-1} \binom{2m-1}{i}=2^{2m-2}.
$$
\end{lemma}

\begin{proof}
We use $\binom{n}{\ell}= \binom{n}{n-\ell}$. If $x:= \sum_{i=0}^m\binom{2m}{m}$, then
$$2x=\sum_{i=0}^{2m}\binom{2m}{i}+\binom{2m}{m} = 2^{2m} +\binom{2m}{m},$$
which gives the first relation. If $y:= \sum_{i=0}^{m-1}\binom{2m-1}{i}$, then
$$2y= \sum_{i=0}^{2m-1}\binom{2m-1}{i}= 2^{2m-1},$$
which gives the second relation.
\end{proof}

\section{Proofs of Theorems \ref{T1.1} to \ref{T1.2b}}\label{sec4}

\noindent{\em Proof of Theorem} \ref{T1.1}.

As already mentioned, the first part of Theorem \ref{T1.1} follows from  Theorem \ref{T1.2a}, which will be proved below.

To prove the second part of Theorem \ref{T1.1}, we assume that $d/N\to \delta$ as $d\to\infty$, where $0\le\delta<1/2$. We write (\ref{3.3}) in the form
\begin{equation}\label{4.12}
\frac{\bE f_k(C_N)}{\binom{N}{k}} = 2^k\frac{\binom{N-k-1}{d-k-1}\left[1+\binom{N-k-1}{d-k-1}^{-1}\sum_{i=0}^{d-k-2}\binom{N-k-1}{i}\right]} {\binom{N-1}{d-1}\left[1+\binom{N-1}{d-1}^{-1}\sum_{i=0}^{d-2}\binom{N-1}{i}\right]},
\end{equation}
and here
$$
\frac{\binom{N-k-1}{d-k-1}}{\binom{N-1}{d-1}}= \frac{d-1}{N-1}\cdots\frac{d-k}{N-k}\to \delta^k\quad\mbox{as }d\to\infty.
$$
Since also $(d-k)/(N-k)\to\delta$, we deduce from Lemma \ref{L4.4} that the normalized sums in the numerator and denominator of (\ref{4.12}) tend to the same finite limit. It follows that $\lim_{d\to\infty} \bE f_k(C_N)/\binom{N}{k}=(2\delta)^k$. \hfill$\Box$

\vspace{2mm}

\noindent{\em Proof of Theorem} \ref{T1.2a}.

As in Section \ref{sec3}, we denote by $\xi_n$ a random variable which has the binomial distribution with parameters $n\in\N$ and $p=\frac{1}{2}$. Let $\xi_n^*$ denote the standardized version of $\xi_n$, that is,
$\xi_n^*=(\xi_n-\bE(\xi_n))/\sqrt{\bV(\xi_n)}$ with $\bE(\xi_n)=\frac{n}{2}$ and $\bV(\xi_n)=\frac{n}{4}$. Then (\ref{2.1b})
implies that
\begin{align*}
\frac{\bE f_k(C_N)}{\binom{N}{k}}&= \frac{\bP(\xi_{N-k-1}\le d-k-1)}{\bP(\xi_{N-1}\le d-1)}
=\frac{\bP\left(\xi_{N-k-1}^*\le \frac{2d-N-k-1}{\sqrt{N-k-1}}\right)}{\bP\left(\xi_{N-1}^*\le \frac{2d-N-1}{\sqrt{N-1}}\right)}\\
&=\frac{\Phi\left(\frac{2d-N-k-1}{\sqrt{N-k-1}}\right)+O\left(\frac{1}{\sqrt{N-k-1}}\right)}{
\Phi\left(\frac{2d-N-1}{\sqrt{N-1}}\right)+O\left(\frac{1}{\sqrt{N-1}}\right)},
\end{align*}
by the Berry--Esseen Theorem (see, e.g., Shiryaev \cite[p. 426]{Shi06}), where $\Phi$ is the distribution function of the standard normal distribution.

We have
$$
\frac{2d-N-k-1}{\sqrt{N-k-1}}=2\sqrt{N}\left(\frac{d}{N}-\frac{1}{2}\right)\frac{1}{\sqrt{1-\frac{k+1}{N}}}-
\frac{k+1}{\sqrt{N-k-1}}.
$$
Since
$$
1-x\le \sqrt{1-x}\le 1-\frac{x}{2} \quad \mbox{for }x\in [0,1],
$$
we get for $\frac{k+1}{N}\le \frac{1}{2}$ (which holds if $d$ is large enough) that
$$
1+\frac{k+1}{2N}\le \frac{1}{\sqrt{1-\frac{k+1}{N}}}\le 1+2\frac{k+1}{N},
$$
thus
$$ \frac{1}{\sqrt{1-\frac{k+1}{N}}} =1+\theta \frac{k+1}{2N}$$
with $\theta\in[2,4]$.
This shows that
\begin{align*}
\frac{2d-N-k-1}{\sqrt{N-k-1}}&=2\sqrt{N}\left(\frac{d}{N}-\frac{1}{2}\right)+\theta\frac{k+1}{{2N}}
2\sqrt{N}\left(\frac{d}{N}-\frac{1}{2}\right)-\frac{k+1}{\sqrt{N-k-1}}\\
&=2\sqrt{N}\left(\frac{d}{N}-\frac{1}{2}\right)+\theta \frac{k+1}{\sqrt{N}}\left(\frac{d}{N}-\frac{1}{2}\right)-
\frac{k+1}{\sqrt{N-k-1}}.
\end{align*}
We define
$$
a(d):=2\sqrt{N}\left(\frac{d}{N}-\frac{1}{2}\right),\qquad b(d):=\theta \frac{k+1}{\sqrt{N}}\left(\frac{d}{N}-\frac{1}{2}\right)-\frac{k+1}{\sqrt{N-k-1}},
$$
hence $a(d)\ge 2a$, $b(d)\to 0$ as $d\to\infty$ and
$$
\frac{2d-N-k-1}{\sqrt{N-k-1}}=a(d)+b(d).
$$
In the same way we get
$$
\frac{2d-N-1}{\sqrt{N-1}}=a(d)+c(d)
$$
with
$$c(d) =\overline\theta \frac{1}{\sqrt{N}}\left(\frac{d}{N}-\frac{1}{2}\right)-\frac{1}{\sqrt{N-1}}  \to 0\quad \mbox{as } d\to\infty,$$
where $\overline{\theta}\in [2,4]$.
Thus we arrive at
\begin{align}\label{4.12a}
\frac{\bE f_k(C_N)}{\binom{N}{k}}&=\frac{\Phi(a(d)+b(d))+O\left(\frac{1}{\sqrt{N-k-1}}\right)}{
\Phi(a(d)+c(d))+O\left(\frac{1}{\sqrt{N-1}}\right)}\nonumber\\
&=\frac{\Phi(a(d))+b(d)\cdot \Phi'(z_1)+O\left(\frac{1}{\sqrt{N-k-1}}\right)}
{\Phi(a(d))+c(d)\cdot\Phi'(z_2)+O\left(\frac{1}{\sqrt{N-1}}\right)}
\end{align}
with intermediate values $z_1,z_2\in\R$. Since the derivative $\Phi'$ is bounded, further $b(d),c(d)\to 0$ as $d\to\infty$, and $\Phi(a(d))\ge \Phi(2a)>0$, we conclude that the quotient tends to $1$ as $d\to\infty$.

For the remaining assertion, we assume that $k(d)$ is bounded. Then we have $b(d)=O(1/\sqrt{N})$. Since the case where $a(d)$ is bounded from below has been settled above, we can assume that $a(d)\le -1$. Thus,
$$
1\le -a(d)=2\sqrt{N}\left(\frac{1}{2}-\frac{d}{N}\right)\le 2a\left(\log\log N\right)^\alpha.
$$
In view of (\ref{4.12a}), it is sufficient to show that
$$
N^{-\frac{1}{2}}\cdot \Phi(a(d))^{-1}=o(1).
$$
To verify this, we use that
$$
1-\Phi(x)\ge \frac{1}{2}\, \frac{1}{x}\, \Phi'(x),\qquad x\ge 1,
$$
(the difference of the function on the left-hand side and the right-hand side converges to zero, as $x\to\infty$, the derivative of this difference is non-positive for $x\ge 1$). Applying this inequality with $x=-a(d)$, we get
\begin{align*}
\frac{1}{\sqrt{N}}\Phi(a(d))^{-1}&=\frac{1}{\sqrt{N}}\frac{1}{1-\Phi(-a(d))}\\
&\le 2\,(-a(d))\,\frac{1}{\sqrt{N}}\,\sqrt{2\pi}\,\exp\left(\frac{1}{2}a(d)^2\right)\\
&=2\sqrt{2\pi}\, \exp\left(\frac{1}{2}a(d)^2+\log(-a(d))-\frac{1}{2}\log N\right).
\end{align*}
Since
\begin{align*}
&\frac{1}{2}a(d)^2+\log(-a(d))-\frac{1}{2}\log N \\
&\qquad
\le 2 a^2(\log\log N)^{2\alpha}+\log(2a)+\alpha\log\log\log N-\frac{1}{2}\log N\to -\infty
\end{align*}
as $d\to\infty$, the assertion follows.
\hfill$\Box$

\vspace{2mm}

\noindent{\em Proof of Theorem} \ref{T1.2b}.

Suppose that $N=2d$. Then by (\ref{3.3}) we have
$$
\frac{\bE f_k(C_N)}{\binom{N}{k}}=2^k\cdot \frac{\sum_{i=0}^{d-k-1}\binom{2d-k-1}{i}}{\sum_{i=0}^{d-1}\binom{2d-1}{i}}.
$$
We distinguish two cases. If $k$ is odd, say $k=2\ell-1$ with $\ell\in\N$, then
\begin{eqnarray*}
\sum_{i=0}^{d-k-1}\binom{2d-k-1}{i} &=& \sum_{i=0}^{d-\ell} \binom{2d-2\ell}{i}-\sum_{i=d-2\ell+1}^{d-\ell}\binom{2d-2\ell}{i}\\
&=& 2^{2d-2\ell-1} + \frac{1}{2}\binom{2d-2\ell}{d-\ell}-\sum_{j=0}^{\ell-1}\binom{2d-2\ell}{d-1-j},
\end{eqnarray*}
where Lemma \ref{L4.2} was used. Since Stirling's formula \eqref{Stirling} yields
$$
 \binom{2d-2\ell}{d-1-j}\sim\frac{4^{d-\ell}}{\sqrt{\pi d}},\enspace j=0,\dots\ell -1,
$$
it follows (again using Lemma \ref{L4.2}) that
\begin{align*}
\frac{\bE f_k(C_N)}{\binom{N}{k}}&\sim 2^{2\ell-1}\frac{1}{4^{d-1}}\left(\frac{1}{2} 4^{d-\ell}+\frac{1}{2}\frac{4^{d-\ell}}{\sqrt{\pi d}}-\ell\frac{4^{d-\ell}}{\sqrt{\pi d}}\right)\\
&= 1-2\cdot 4^{\ell-d}\left(\ell-\frac{1}{2}\right)\frac{4^{d-\ell}}{\sqrt{\pi d}}=1-(2\ell-1)\frac{1}{\sqrt{\pi d}}.
\end{align*}
If $k$ is even, say $k=2\ell$  with $\ell\in\N$, then
\begin{eqnarray*}
\sum_{i=0}^{d-k-1}\binom{2d-k-1}{i} &=& \sum_{i=0}^{d-\ell-1}\binom{2(d-\ell)-1}{i}-\sum_{j=1}^\ell\binom{2(d-\ell)-1}{d-j}\\
&=& 2^{2d-2\ell-2} - \sum_{j=1}^\ell\binom{2(d-\ell)-1}{d-j}
\end{eqnarray*}
by Lemma \ref{L4.2}. By Stirling's approximation \eqref{Stirling},
$$
\binom{2(d-\ell)-1}{d-j}\sim \frac{1}{2}\frac{4^{d-\ell}}{\sqrt{\pi d}},\enspace j=1,\ldots,\ell,
$$
and hence we get
\begin{eqnarray*}
\frac{\bE f_k(C_N)}{\binom{N}{k}}
&\sim& 2^{2\ell}\frac{1}{4^{d-1}}\left(2^{2d-2\ell-2} -\ell\frac{1}{2}\frac{4^{d-\ell}}{\sqrt{\pi d}}\right)\\
&=& 1-(2\ell)\frac{1}{\sqrt{\pi d}}.
\end{eqnarray*}
Thus in both cases the asymptotic relation is proved. \hfill$\Box$

\section{Proofs of Theorems \ref{T1.3} to \ref{T1.4b}}\label{sec5}

\noindent{\em Proof of Theorem} \ref{T1.3}.

The first part of Theorem \ref{T1.3} follows from Theorem \ref{T1.4a}, which will be proved below.

For the second part of the proof, we assume that $d/N\to \delta$ as $d\to\infty$, with $0\le\delta<1/2$. We note that relation (\ref{3.3a}) shows that
$$
\bE 2 U_{d-k}(C_N)=1- \frac{C(N,d-k)}{C(N,d)}.
$$
Let
$$
R(N-1,d-\ell):=\frac{1}{\binom{N-1}{d-\ell}}\sum_{i=0}^{d-\ell-1}\binom{N-1}{i},\quad \ell\in\{0,\ldots,d-1\}.
$$
Then
$$
\frac{C(N,d-k)}{C(N,d)}=\frac{\binom{N-1}{d-k}}{\binom{N-1}{d}}\frac{R(N-1,d-k)}{R(N-1,d)},
$$
and here
$$
\frac{\binom{N-1}{d-k}}{\binom{N-1}{d}}=\frac{d}{N-d}\cdots\frac{d-(k-1)}{N-d+k-1} \to \left(\frac{\delta}{1-\delta}\right)^k\quad\mbox{as }d\to\infty.
$$

Together with $d/N\to\delta$ we have $(d-k+1)/N\to\delta$ and $(d+1)/N\to\delta$, hence Lemma \ref{L4.4} yields that
$$
R(N-1,d-k)\to\frac{\delta}{1-2\delta}\quad\text{and}\quad R(N-1,d)\to\frac{\delta}{1-2\delta}.
$$
The assertion follows.  \hfill$\Box$

\vspace{2mm}

\noindent{\em Proof of Theorem} \ref{T1.4a}.

The random variables $\xi_n$ and $\xi_n^*$ are defined as in the proof of Theorem \ref{T1.2a}. Using (\ref{2.1c}), we proceed as in that proof and obtain
\begin{align*}
1-\bE 2U_{d-k}(C_N)&= \frac{\bP(\xi_{N-1}\le d-k-1)}{\bP(\xi_{N-1}\le d-1)}
=\frac{\bP\left(\xi_{N-1}^*\le \frac{2d-N-k-1}{\sqrt{N-k-1}}\right)}{\bP\left(\xi_{N-1}^*\le \frac{2d-N-1}{\sqrt{N-1}}\right)}\\
&=\frac{\Phi\left(\frac{2d-N-k-1}{\sqrt{N-k-1}}\right)+O\left(\frac{1}{\sqrt{N-1}}\right)}{
\Phi\left(\frac{2d-N-1}{\sqrt{N-1}}\right)+O\left(\frac{1}{\sqrt{N-1}}\right)}.
\end{align*}
The rest of the proof follows that for Theorem \ref{T1.2a}. \hfill$\Box$

\vspace{2mm}

\noindent{\em Proof of Theorem} \ref{T1.4b}.

Let $N=2d$. Using Lemma \ref{L4.2}, we get
$$
\frac{C(2d,d-k)}{C(2d,d)}=\frac{\sum_{i=0}^{d-1}\binom{2d-1}{i}-\sum_{i=d-k}^{d-1}\binom{2d-1}{i}}{
\sum_{i=0}^{d-1}\binom{2d-1}{i}}=1-4^{1-d}\sum_{j=0}^{k-1}\binom{2d-1}{d-1-j}.
$$
By Stirling's approximation \eqref{Stirling}, we obtain for $j\in\{0,\ldots,k-1\}$ that
$$
\binom{2d-1}{d-1-j}
\sim \frac{4^d}{2}\frac{1}{\sqrt{\pi d}}.
$$
Hence we get
$$
\frac{C(2d,d-k)}{C(2d,d)}\sim 1-4^{1-d}\cdot k\cdot \frac{1}{2}4^d\frac{1}{\sqrt{\pi d}}=1-2k\frac{1}{\sqrt{\pi d}}.
$$
Thus we arrive at
$$
\bE U_{d-k}(C_{2d})\sim \frac{1}{2}-\frac{1}{2}\left(1-2k\frac{1}{\sqrt{\pi d}}\right)=\frac{k}{\sqrt{\pi d}},
$$
which completes the proof. \hfill$\Box$

\section{Linearly growing face dimensions}\label{sec6}

In this section and the next, we allow also the number $k$ to grow linearly with the dimension $d$. In the present section, we are interested in a phase transition for the expectation $\bE f_k(C_N)/\binom{N}{k}$. It turns out that it appears at the same threshold as it was observed earlier by Donoho and Tanner \cite{DT10} for a different, but closely related class of random polyhedral cones. These authors considered a real random $d\times N$ matrix ${\sf A}$ of rank $d$, where $d<N$, the nonnegative orthant
$$ \R^N_+:=\{x\in\R^N:x_i\ge 0 \mbox{ for }i=1,\dots,N\}$$
of $\R^N$, and its image ${\sf A}\R^N_+$ in $\R^d$. Considering the column vectors of ${\sf A}$ as random vectors in $\R^d$, the image ${\sf A}\R^N_+$ is the positive hull of these vectors. For a suitable distribution, the random cone ${\sf A}\R^N_+$ is obtained in a similar way as the Cover--Efron cone, just by omitting the condition that the cone is different from $\R^d$. Imposing this condition leads, of course, to different distributions of the random cones. Comparing formula (13) of \cite{DT10} with our formula (\ref{3.3}), where the right-hand side can be written as
$$ \frac{1-P_{N-d,N-k}}{P_{d,N}},$$
we see that it results in an additional denominator in the expression for the expected number of $k$-faces, thus increasing this expectation. Therefore some of our estimates, though leading to the same threshold, require more effort.

\vspace{2mm}

\noindent{\em Proof of Theorem} \ref{T1.5}.

First we assume that $0\le \rho<\rho_W(\delta)$. For this part of the proof, we reproduce an argument which was suggested by an anonymous referee (our original proof can be found in arXiv:2004.11473v1).

We use the representation
\begin{equation}\label{neweq1}
\frac{\bE f_k(C_N)}{\binom{N}{k}} =\frac{\bP(\xi_{N-k-1}\le d-k-1)}{\bP(\xi_{N-1}\le d-1)},
\end{equation}
given by (\ref{2.1b}), where $\xi_n$ is a random variable with binomial distribution with parameters $n$ and $1/2$. Since $\xi_n$ is in distribution equal to the sum of $n$ i.i.d. Bernoulli random variables with parameter $1/2$, the weak law of large numbers gives
$$ \lim_{n\to\infty} \bP\left(\xi_n>n\left(\frac{1}{2}+\varepsilon\right)\right)=0$$
for each $\varepsilon>0$.

If now $\frac{d}{N}\to\delta>\frac{1}{2}$ and $\frac{k}{d}\to\rho<2-\delta^{-1}$, then
$$ \lim_{d\to\infty} \frac{d-k-1}{N-k-1}= \frac{1-\rho}{1/\delta-\rho} >\frac{1}{2} \quad\mbox{and}\quad \lim_{d\to\infty}\frac{d-1}{N-1} =\delta>\frac{1}{2}.$$
Therefore,
$$ \lim_{d\to\infty} \bP(\xi_{N-k-1}\le d-k-1)=1,\qquad \lim_{d\to\infty} \bP(\xi_{N-1}\le d-1)=1$$
and hence
$$ \lim_{d\to\infty} \frac{\bE f_k(C_N)}{\binom{N}{k}} =1,$$
as stated.

For the second part of the proof, we use (\ref{4.12}) and show first that for increasing $d$ the terms in brackets remain between two positive constants. The asymptotic behavior of the remaining quotient is then determined with the aid of Stirling's formula.

We assume that $\rho>\rho_W(\delta)$. Then $\rho>2-\delta^{-1}$, irrespective of whether $\delta\ge 1/2$ or not.
For sufficiently large $d$ (which we assume in the following), we then have $\rho_d>2-\delta_d^{-1}$, which implies $N-2d+k>0$. Thus, we can apply ({\ref{4.2}) to the normalized sum in the numerator of (\ref{4.12}). This yields
$$ \frac{d-k-1}{N-d+1} \le \binom{N-k-1}{d-k-1}^{-1}\sum_{i=0}^{d-k-2}\binom{N-k-1}{i} \le \frac{d-k-1}{N-d+1} \cdot\frac{N-d+2}{N-2d+k+4}.$$
Here,
$$ \lim_{d\to\infty} \frac{d-k-1}{N-d+1} = \frac{\delta(1-\rho)}{1-\delta},\qquad  \lim_{d\to\infty} \frac{N-d+2}{N-2d+k+4} = \frac{1-\delta}{1-2\delta+\rho\delta},$$
where the last denominator is positive. It follows that
\begin{equation}\label{Eqcc}
 c_1\le \binom{N-k-1}{d-k-1}^{-1}\sum_{i=0}^{d-k-2}\binom{N-k-1}{i} \le c_2
\end{equation}
for all sufficiently large $d$. Here and below we denote by $c_i$ a positive constant that is independent of $d$.

In view of (\ref{4.12}), we now determine the asymptotic behavior of
$$ 2^k\frac{\binom{N-k-1}{d-k-1}}{\binom{N-1}{d-1}} = 2^k\frac{N}{d}\frac{d-k}{N-k}\frac{\binom{d}{k}}{\binom{N}{k}} \quad\mbox{as }d\to\infty.$$
Here,
$$ \lim_{d\to\infty} \frac{N}{d}\frac{d-k}{N-k} =\frac{1-\rho}{1-\rho\delta}.$$
To treat the remaining terms, we use the Stirling formula \eqref{Stirling}. Define

$$ \frac{d}{N}=: \delta_d,\qquad \frac{k}{d}=:\rho_d,$$
then $\delta_d\to\delta$ and $\rho_d\to\rho$ as $d\to\infty$. We obtain
\begin{align}
 2^k\frac{\binom{d}{k}}{\binom{N}{k}}
 &= 2^{\rho_d\delta_d N} \sqrt{\frac{1-\rho_d\delta_d}{1-\rho_d}}\cdot\frac{(\delta_d N)^{\delta_d N} (N-\rho_d \delta_d N)^{N-\rho_d\delta_d N}}{(\delta_d N -\rho_d\delta_d N)^{\delta_d N-\rho_d\delta_d N}N^N} \cdot e^{\varphi/12 N}\nonumber\\
 &=   e^{\varphi/12 N}\sqrt{\frac{1-\rho_d\delta_d}{1-\rho_d}}H(\delta_d,\rho_d)^N ,\label{eq01}
 \end{align}
where $\varphi$ is contained in a fixed interval independent of $d$, and
$$
H(a,b):=\frac{(2a)^{ab}(1-ab)^{1-ab}}{(1-b)^{a(1-b)}},\quad a\in (0,1),\,b<1.
$$
We define
\begin{equation}\label{g}
g(a):= H\left(a,2-a^{-1}\right)= 2a^a(1-a)^{1-a} \qquad\mbox{for }a\in (0,1).
\end{equation}
Note that for $a,b\in (0,1)$ we have $b>2-a^{-1}$ if and only if $a<1/(2-b)$.
  Let $H_a(b):=H(a,b)$. Differentiation yields
$$ H_a'(b)=a\log\left(\frac{2a(1-b)}{1-ab}\right) H(a,b).$$
Hence $H_a'(b)<0$ for $b>2-a^{-1}$, since
$$
\frac{2a(1-b)}{1-ab}<1 \Leftrightarrow b>2-a^{-1}.
$$
If $a\le 1/2$, then
$$ H(a,b)<H(a,0)=1$$
for $b\in (0,1)$. On the other hand, if $a>1/2$ and $b>2-a^{-1}$, we  have
\begin{equation}\label{6.x}
H(a,b)< H\left(a,2-a^{-1}\right)=2a^a(1-a)^{1-a}.
\end{equation}
Since the function $g$ defined by (\ref{g}) satisfies $g(1/2)=1$ and $g'(a)>0$ for $1/2<a<1$, we also have
\begin{equation}\label{6.y}
H(a,2-a^{-1})>1\quad \mbox{if }a>1/2.
\end{equation}
Now we distinguish two cases.

(1) Let $\delta\le 1/2$. Then \eqref{4.12}, \eqref{Eqcc} and \eqref{eq01} yield
\begin{align*}
\limsup_{d\to\infty}\frac{\bE f_k(C_N)}{\binom{N}{k}}&\le
\frac{1-\rho}{1-\rho\delta}\frac{\sqrt{1-\rho\delta}}{\sqrt{1-\rho}}
\limsup_{d\to\infty} H(\delta_d,\rho_d)^N\frac{1+c_2}{1}\\
&\le (1+c_2)\limsup_{d\to\infty}c_3^N=0,
\end{align*}
where $H(\delta_d,\rho_d)\le c_3<1$, since $H(\delta_d,\rho_d)\to H(\delta,\rho)<H(\delta,0)=1$.

(2) Let $\delta> 1/2$. Then we can assume that $N/d<c_4<2$. We have
$$
\sum_{i=0}^{d-2}\binom{N-1}{i}=2^{N-1}-\sum_{j=0}^{N-d}\binom{N-1}{j}.
$$
Since $2(N-d)<N$,  Lemma \ref{L4.1} yields
$$
\binom{N-1}{N-d+1}^{-1}\sum_{j=0}^{N-d}\binom{N-1}{j}\le \frac{\frac{N}{d-1}-1}{2-\frac{N}{d}},
$$
and hence
$$
\binom{N-1}{d-1}^{-1}\sum_{i=0}^{d-2}\binom{N-1}{i}\ge\binom{N-1}{d-1}^{-1} 2^{N-1}-\frac{1}{2-\frac{N}{d}}.
$$
To estimate the last binomial coefficient, we use Stirling's approximation \eqref{Stirling} together with (\ref{6.x}). Thus, we get for large $d$ the lower bound
$$
\binom{N-1}{d-1}^{-1}\sum_{i=0}^{d-2}\binom{N-1}{i}\ge c_5 H\left(\delta_d,2-\delta_d^{-1}\right)^N-c_6.
$$
Combining these estimates and starting again from \eqref{4.12}, we finally obtain
\begin{align*}
\limsup_{d\to\infty}\frac{\bE f_k(C_N)}{\binom{N}{k}}&\le c_7 \limsup_{d\to\infty} H(\delta_d,\rho_d)^N\frac{1+c_2}{
1+ c_5H(\delta_d,2-\delta_d^{-1})^N -c_6}\\
&= c_8 \limsup_{d\to\infty}\left(\frac{H(\delta_d,\rho_d)}{H(\delta_d,2-\delta_d^{-1})}\right)^N\frac{1}{(1-c_6)H(\delta_d,2-\delta_d^{-1})^{-N}+c_5}\\
&=0.
\end{align*}
Here we have used that
$$
\frac{H(\delta,\rho)}{H(\delta,2-\delta^{-1})}<1\quad \text{for }
\rho>2-\delta^{-1}
$$
by (\ref{6.x}) and that
$$
H(\delta,2-\delta^{-1})>1 \quad \mbox{for }\delta > \frac{1}{2}$$
by (\ref{6.y}). This completes the proof also in the case $\rho>\rho_W(\delta)$.
\hfill$\Box$

\vspace{3mm}

\noindent
\textbf{Remark.} Under the assumption $\frac{d}{N}\to\delta>\frac{1}{2}$ as $d\to\infty$, we have seen in the first part of the preceding proof that $
\bP\left({\rm pos}\{X_1,\dots,X_{N-1}\}\not=\R^d\right)=\bP(\xi_{N-1}\le d-1)\to 1$ for $d\to\infty$,
 as a simple consequence of the weak law of large numbers. In this situation, an application of a large deviation (concentration) result for the binomial distribution in fact shows that the convergence is exponentially fast. For this, we choose $d$ sufficiently large so that $\frac{d}{N-1}\ge \frac{1}{2}$. Then
\begin{align*}
\bP(\xi_{N-1}\ge d)&=\bP\left(\xi_{N-1}-\frac{N-1}{2}\ge\left(\frac{d}{N-1}- \frac{1}{2}\right)(N-1)\right)\\
&\le \exp\left(-2\left(\frac{d}{N-1}- \frac{1}{2}\right)^2(N-1)\right),
\end{align*}
by Okamoto's inequality (see \cite[Theorem 2 (i)]{Okamoto}), which applies since $\frac{d}{N-1}- \frac{1}{2}\ge 0$ (if $d$ is large enough). Hence, if $d$ is sufficiently large so that
$\frac{d}{N-1}- \frac{1}{2}\ge\frac{1}{\sqrt{2}}\left(\delta-\frac{1}{2}\right)>0$, we obtain
$$
\bP(\xi_{N-1}\ge d)\le \exp\left(-\left(\delta-\frac{1}{2}\right)^2(N-1)\right),
$$
as asserted. On the other hand, if $\frac{d-1}{N-1}\le \frac{1}{2}$ then \cite[Theorem 2 (ii)]{Okamoto} yields
\begin{align*}
\bP(\xi_{N-1}\le d-1)&=\bP\left(\xi_{N-1}-\frac{N-1}{2}\le-\left( \frac{1}{2}-\frac{d-1}{N-1}\right)(N-1)\right)\\
&\le \exp\left(-2\left(\frac{d}{N-1}- \frac{1}{2}\right)^2(N-1)\right).
\end{align*}
Hence, if $\frac{d}{N}\to\delta<\frac{1}{2}$ and $d$ is large enough so that $\frac{1}{2}-\frac{d}{N-1} \ge\frac{1}{\sqrt{2}}\left(\frac{1}{2}-\delta\right)>0$, then
$$
\bP(\xi_{N-1}\le d-1)\le \exp\left(-\left(\delta-\frac{1}{2}\right)^2(N-1)\right).
$$
In this case, a finer analysis of the  ratio \eqref{neweq1}  with corresponding lower bounds for the involved probabilities is required.

\vspace{3mm}

We prepare the proof of Theorem \ref{T1.6} by a lemma, which serves to establish the threshold $\rho_S$ and to provide an upper estimate for it. By ${\sf H}$ we denote the binary entropy function with base $e$, that is
$$ {\sf H}(x):= -x\log x-(1-x)\log(1-x) \quad\mbox{for }0\le x\le 1$$
(with $0\log 0:=0$). We note that ${\sf H}(0)={\sf H}(1)=0$ and that ${\sf H}$ attains its unique maximum, $\log 2$, at the point $1/2$.
As in \cite{DT10}, we consider the function defined by
$$
G(\delta,\rho):={\sf H}(\delta)+\delta {\sf H}(\rho)-(1-\rho\delta )\log 2,\quad \rho,\delta\in[0,1].
$$
For a later application, we remark that
\begin{equation}\label{6.z}
e^{-G(\delta,\rho)} = (1-\delta)^{1-\delta}\delta^{\delta}(1-\rho)^{\delta(1-\rho)}\rho^{\delta\rho}2^{1-\delta\rho}.
\end{equation}

\begin{lemma}\label{L6.1}
For $\delta\in (1/2,1)$, the function $G_\delta$ defined by $G_\delta(x):= G(\delta,x)$ has a unique zero $x_0\in [0,1]$. Moreover, $x_0\in (0,\min\{\frac{2}{3},2-\delta^{-1}\})$.
\end{lemma}

\begin{proof}
Clearly, $G_\delta(0)={\sf H}(\delta)-\log 2<0$ since $\delta\not=1/2$. We have
$$
G_\delta'(x)=\delta \log\left(\frac{2(1-x)}{x}\right).
$$
Hence $x_0=2/3$ is the unique zero of $G_\delta'$ in $(0,1)$, and $G_\delta'>0$ in $(0,2/3)$ and $G_\delta'<0$ in $(2/3,1)$. We will show that

(a) $G_\delta (2/3)>0$,

(b) $G_\delta(1)>0$,

(c) $G_\delta (2-\delta^{-1})>0$,

\noindent
which then implies that $G_\delta$ has a unique zero $x_0$ in $[0,1]$ and $x_0<2/3$ as well as  $x_0<2-\delta^{-1}$.

For (a) we define $\kuppa_1(\delta):=G_\delta(2/3)$ for $\delta\in [1/2,1]$. Then a simple calculation shows that
$\kuppa_1(1/2)=1/2\log 3>0$, $\kuppa_1(1)=\log3-\log 2>0$ and
$$
\kuppa_1'(\delta)=\log\left(\frac{3(1-\delta)}{\delta}\right),\quad \delta\in (1/2,1).
$$
Hence $\kuppa_1'(3/4)=0$, $\kuppa_1'>0$ on $(1/2,3/4)$ and $\kuppa_1'<0$ on $(3/4,1)$. This shows that $\kuppa_1(\delta)=G_\delta (2/3)>0$ for $\delta\in (1/2,1)$ (with maximal value $\kuppa_1(3/4)=\log 2$).

For (b)  we consider $\kuppa_2(\delta):=G_\delta(1)$ for $\delta\in [1/2,1]$.
 We have $\kuppa_2(1/2)=1/2\log 2>0$, $\kuppa_2(1)=0$ and
$$
\kuppa_2'(\delta)=\log\left(\frac{2(1-\delta)}{\delta}\right),\quad \delta\in (1/2,1).
$$
Hence, $\kuppa_2'(2/3)=0$, $\kuppa_2'>0$ on $(1/2,2/3)$ and $\kuppa_2'<0$ on $(2/3,1)$. In particular, this yields $G_\delta(1)=\kuppa_2(\delta)>0$ for $\delta\in (1/2,1)$.

For (c) we consider $\kuppa_3(\delta)=G_\delta(2-\delta^{-1})$ for $\delta\in [1/2,1]$. Then
$$
\kuppa_3(\delta)=-2(1-\delta)\log(1-\delta)-(2\delta-1)\log(2\delta-1)-2(1-\delta)\log 2
$$
and
$$
\kuppa_3'(\delta)=2\log\left(\frac{2(1-\delta)}{2\delta-1}\right),\quad \delta\in (1/2,1).
$$
We have $\kuppa_3(1/2)=\kuppa_3(1)=0$, $\kuppa_3'(3/4)=0$, $\kuppa_3'>0$ on $(1/2,3/4)$ and $\kuppa_3'<0$ on $(3/4,1)$. But then
$\kuppa_3>0$ on $(1/2,1)$ and therefore
$G_\delta(2-\delta^{-1})=\kuppa_3(\delta)>0$ for $\delta\in (1/2,1)$.
\end{proof}

Now we denote the zero $x_0$ of the function $G_\delta$ provided by Lemma \ref{L6.1} by $\rho_S(\delta)$. This defines the function $\rho_S$ appearing in Theorem \ref{T1.6}.

\vspace{3mm}

\noindent{\em Proof of Theorem} \ref{T1.6}.

Again we define
$$ \frac{d}{N}=: \delta_d,\qquad \frac{k}{d}=:\rho_d.$$
Let $\delta>1/2$ and $0<\rho<\rho_S(\delta)$. Then $G(\delta,\rho)<0$. For sufficiently large $d$, we have $G(\delta_d,\rho_d)<0$ as well as $\delta_d>1/2$ and $\rho_d<\rho_S(\delta_d)$. We assume that $d$ is large enough in this sense. Since $\rho_s(\delta_d)<\rho_W(\delta_d)$, we have $N-2d+k<0$. We can, therefore, apply the estimates from the first part of the proof of Theorem \ref{T1.5}.

Let $X_1,\dots,X_N$ be i.i.d.~unit vectors with distribution $\phi$. By definition, $C_N$ is the positive hull of $X_1,\dots,X_N$ under the condition that this positive hull is different from $\R^d$. For $k\in \{1,\dots,d-1\}$, choose $1\le i_1<\dots<i_k\le N$ and let $M=\{X_{i_1},\dots,X_{i_k}\}$. Then the distribution of ${\rm pos}\,M$, under the condition that ${\rm pos}\{X_1,\dots,X_N\}\not=\R^d$, is independent of the choice of $i_1,\dots,i_k$, hence
$$ \binom{N}{k}\bP\left({\rm pos}\,M\in \F_k(C_N)\right)= \bE f_k(C_N).$$
Therefore,
$$ p:= \bP\left({\rm pos}\,M\notin \F_k(C_N)\right)= 1-\frac{\bE f_k(C_N)}{\binom{N}{k}} = \frac{A}{1+A}$$
by (\ref{3.1a}) (and with the notation used there). By Boole's inequality,
$$ \bP\left({\rm pos}\,M\notin \F_k(C_N) \mbox{ for some $k$-element subset }M\right) \le \binom{N}{k}p$$
and thus
$$ \bP\left(f_k(C_N)=\binom{N}{k}\right) \ge 1-\binom{N}{k}p.$$
Here,
$$ \binom{N}{k}-\bE f_k(C_N)= \binom{N}{k}p\le \binom{N}{k} A.$$
By Lemmas \ref{L3.1} and \ref{L3.2},
\begin{equation}\label{n3}
A^{-1} = \frac{2^{N-1}P_{d-k,N-k}}{\sum_{j=1}^k\binom{k}{j}\sum_{m=0}^{j-1}\binom{N-k-1}{d-k+m}}
\ge \frac{2^{N-2}-2^{k-1}\sum_{r=0}^{N-2d+k-1} \binom{N-k-1}{d-k+r}} {\sum_{j=1}^k\binom{k}{j}\sum_{m=0}^{j-1}\binom{N-k-1}{d-k+m}}.
\end{equation}
Since $N-2d+k-1< 0$ for all $d$, the sum in the last numerator of (\ref{n3}) is zero, hence
\begin{equation}\label{24}
A\le {\sum_{j=1}^k\binom{k}{j}\sum_{m=0}^{j-1}2^{2-N}\binom{N-k-1}{d-k+m}}.
\end{equation}
Using this and the identities
\begin{equation}\label{24b}
\binom{N-k-1}{d-k}=\frac{N-d}{N-k}\cdot \frac{(N-k)!}{(N-d)!(d-k)!},\qquad \sum_{j=1}^k\binom{k}{j}j=2^{k-1}\cdot k
\end{equation}
together with \eqref{Stirling}, we get
\begin{align*}
\binom{N}{k}A&\le 2^{2-N}\frac{N-d}{N-k}k2^{k-1}\frac{N!}{(N-d)!(d-k)!k!}\\
&= 2\frac{(1-\delta_d)\delta_d\rho_d}{1-\delta_d\rho_d}N\frac{\sqrt{N}}{2\pi\sqrt{N-d}\sqrt{d-k}\sqrt{k}}\,
e^{G(\delta_d,\rho_d)N}e^\frac{\varphi}{12N},
\end{align*}
where (\ref{6.z}) was used and where $\varphi\in(0,1)$. Since $G(\delta_d,\rho_d)\to G(\delta,\rho) <0$ as $d\to\infty$,  it follows that
$$  \lim_{d\to\infty} \binom{N}{k}A=0,$$
from which the assertions follow. \hfill$\Box$

\section{Proof of Theorem \ref{T1.7}} \label{sec7}

We use the representations
$$
\bE 2 U_{d-k}(C_N)=1- \frac{C(N,d-k)}{C(N,d)}=1-\frac{\bP(\xi_{N-1}\le d-k-1)}{\bP(\xi_{N-1}\le d-1)}
$$
and show the convergence of the quotients, under different assumptions.

To prove the first part of the theorem, we assume that $0\le \rho< \frac{1}{2}\rho_W(\delta)$; then $0\le \rho<1-(2\delta)^{-1}$ and $\delta(1-\rho)>1/2$ and $\delta>1/2$. Hence
$$
\frac{d-k-1}{N-1}\to \delta(1-\rho)>\frac{1}{2}\qquad\text{and}\qquad \frac{d-1}{N-1}\to\delta>\frac{1}{2}.
$$
Using the weak law of large numbers, as in the proof of Theorem \ref{T1.5}, we obtain
$$ \lim_{d\to\infty} \bP(\xi_{N-1}\le d-k-1)=1,\qquad \lim_{d\to\infty} \bP(\xi_{N-1}\le d-1)=1,$$
which completes this part of the argument.

Now we deal with the second part of the proof and point out that our argument requires to distinguish whether $\rho>\rho_W(\delta)$ or not. We begin with the case $\rho>\rho_W(\delta)$; then $\rho\ge 2-\delta^{-1}$. Clearly,
$$
\frac{C(N,d-k)}{C(N,d)}=\frac{\binom{N-1}{d-k-1}}{\binom{N-1}{d-1}}\cdot \frac{1+\binom{N-1}{d-k-1}^{-1}\sum_{i=0}^{d-k-2}
\binom{N-1}{i}}{1+\binom{N-1}{d-1}^{-1}\sum_{i=0}^{d-2}\binom{N-1}{i}}.
$$
Since (a fortiori) $\rho>1-(2\delta)^{-1}$, we have $N-2d+2k>0$ for sufficiently large $d$, hence Lemma \ref{L4.1} yields
\begin{equation}\label{cite1}
\frac{1}{\binom{N-1}{d-k-1}}\sum_{i=0}^{d-k-2}
\binom{N-1}{i}\le \frac{d-k-1}{N-d+k+1}\cdot\frac{N-d+k+2}{N-2d+2k+4}\to \frac{\delta(1-\rho)}{1-2\delta(1-\rho)}
\end{equation}
as $d\to\infty$, and the last denominator is positive. Hence, if $d$ is large enough (which is always assumed in the following), there are constants $C_1,C_2$, independent of $d$, such that
\begin{align*}
0\le \frac{C(N,d-k)}{C(N,d)}&\le C_1\frac{d-k}{d}\frac{d!(N-d)!}{(d-k)!(N-d+k)!}\\
&\le C_2\frac{d^d(N-d)^{N-d}}{(d-k)^{d-k}(N-d+k)^{N-d+k}}\le C_2\cdot K(\delta_d,\rho_d)^N,
\end{align*}
where
$$
K(a,b):=\frac{g(a)}{g(a(1-b))}=\frac{a^{ab}(1-a)^{1-a}}{(1-b)^{a(1-b)}(1-a(1-b))^{1-a(1-b)}},\quad a\in (0,1),\, b\in [0,1).
$$
We have $K(a,0)=1$, and also $K(a,2-a^{-1})=1$ if $a\ge 1/2$. We write $K_a(b):=K(a,b)$. Then
$$
K_a'(b)=a\log\left(\frac{a(1-b)}{1-a(1-b)}\right)K(a,b)<0 \quad \mbox{for}\quad b> 1-(2a)^{-1} \mbox{ (where $b\ge 0$)}.
$$
Thus, for any $a\in (0,1)$, we have $K(a,b)<1$ for all $b>2-a^{-1}$ in $(0,1)$.

Since $\rho>2-\delta^{-1}$, we deduce that $K(\delta,\rho)<1$ and hence that $K(\delta_d,\rho_d)\le c <1$ for all sufficiently large $d$, with $c$ independent of $d$. It follows that
\begin{equation}\label{convzero}
\lim_{d\to\infty}\frac{C(N,d-k)}{C(N,d)} =0.
\end{equation}

Now we suppose that $ \frac{1}{2}\rho_W(\delta) < \rho \le \rho_W(\delta)$, then $1-(2\delta)^{-1}<\rho\le 2-\delta^{-1}$. Since $\rho>0$, we have $\delta>1/2$, further $\delta(1-\rho)<1/2$.

We use repeatedly that
$$ 2^{N-1} = \sum_{i=0}^{d-1}\binom{N-1}{i} + \sum_{i=0}^{N-d-1} \binom{N-1}{i}.$$

We note that still $N-2d+2k>0$ for sufficiently large $d$, so that \eqref{cite1} can be applied. It yields
\begin{eqnarray*}
\frac{C(N,d-k)}{C(N,d)} &=& \frac{\binom{N-1}{d-k-1}\left[1+\binom{N-1}{d-k-1}^{-1} \sum_{i=0}^{d-k-2} \binom{N-1}{i}\right]}{2^{N-1}-\sum_{i=0}^{N-d-1}\binom{N-1}{i}}\\
&\le& \frac{\binom{N-1}{d-k-1}}{2^{N-1}}\cdot\frac{1+C_3}{1-\frac{1}{2^{N-1}}\sum_{i=0}^{N-d-1}\binom{N-1}{i}}.
\end{eqnarray*}
Here and below, $C_m$ denotes a positive constant independent of $d$.

To estimate the last denominator, we can again use Lemma \ref{L4.1}, since $\delta>1/2$ and hence
$2(N-d-1)<N$, if $d$ is large enough, to get
\begin{align*}
\frac{1}{2^{N-1}}\sum_{i=0}^{N-d-1}\binom{N-1}{i} &=\frac{\binom{N-1}{N-d}}{2^{N-1}}\frac{1}{\binom{N-1}{N-d}}\sum_{i=0}^{N-d-1}\binom{N-1}{i}\\
&\le \frac{\binom{N-1}{N-d}}{2^{N-1}}\frac{N-d}{d}\cdot\frac{d+1}{-N+2d+2}\\
&\le C_4\frac{\binom{N-1}{d-1}}{2^{N-1}}\le C_5 g(\delta_d)^{-N}
\end{align*}
with $g$ defined by (\ref{g}). As already observed, $g(1/2)=1$ and $g(a)>1$ for $a\in [0,1]\setminus\{1/2\}$.
Since $\delta>1/2$, we have $g(\delta)>1$ and hence $g(\delta_d)\ge c>1$ for sufficiently large $d$, with $c$ independent of $d$. It follows that $g(\delta_d)^{-N}\to 0$ as $d\to\infty$.

Finally, we observe that
$$
\frac{\binom{N-1}{d-k-1}}{2^{N-1}}\le C_6\frac{N^N}{2^N(d-k)^{d-k}(N-d+k)^{N-d+k}}=C_7\cdot g(\delta_d(1-\rho_d))^{-N}\to 0,
$$
since $\delta_d(1-\rho_d)\to\delta(1-\rho)<1/2$ and  $g(\delta(1-\rho))>1$. Thus, (\ref{convzero}) is obtained again.

\hspace*{\fill}$\Box$

\vspace{2mm}

\noindent{\bf Appendix.} Before Theorem \ref{T1.5}, we have formulated that the occurrence of the same threshold in the work \cite{DT10} and in our Theorem \ref{T1.5} may be unexpected. An anonymous referee suggested the following explanation. We quote it verbally (but adding bibliographic information): ``Here is an attempt of explanation of this coincidence. Donoho and Tanner \cite{DT10} consider projections of orthants on random uniform subspaces. By the same argument as in the paper of Baryshnikov and Vitale \cite{BV94}, the expected number of faces does not change if random uniform projection is replaced by applying a Gaussian random matrix. Since the orthant is the positive hull of the standard basis, it follows that instead of the Donoho-Tanner cones one may consider positive hulls of $N$ i.i.d. standard Gaussian random variables in $\R^d$. Thus, the difference between the Cover-Efron cones studied here and the Donoho-Tanner cones is the conditioning on the event that the cone is not equal to $\R^d$. (\dots) Now let us look at the case $\delta>1/2$ in Theorem \ref{T1.5}. Then, $N$ is between $d$ and $(2-\varepsilon)d$, which means that the probability that the positive hull of the vectors $X_1,\dots,X_N$ is $\R^d$ goes to $0$ exponentially fast. So, both models of cones differ just on an event of exponentially small probability and are equal otherwise. Moreover, the events of exponentially small probability make no contribution to the expected $f$-vector because on this event the Donoho-Tanner cone is $\R^d$ and $f_k=0$ for all $k<d$. The case $\delta <1/2$ is more difficult to explain.`` (end of quotation) The referee then sketches an argument for this case which, in his/her opinion, is not rigorous, but makes the result quite natural.

\vspace{2mm}

\noindent{\bf Statement.} 
Data sharing not applicable to this article as no datasets were generated or analysed during the current study.

\vspace{2mm}

\noindent{\bf Acknowledgment.} We are very grateful to the anonymous referees, whose hints allowed us in some cases to give shorter proofs and to obtain stronger results.

\noindent Authors' addresses:\\[2mm]
Daniel Hug\\Karlsruhe Institute of Technology (KIT), Department of Mathematics\\D-76128 Karlsruhe, Germany\\E-mail: daniel.hug@kit.edu\\[2mm]
Rolf Schneider\\Mathematisches Institut, Albert-Ludwigs-Universit{\"a}t\\D-79104 Freiburg i.~Br., Germany\\E-mail: rolf.schneider@math.uni-freiburg.de


\begin{thebibliography}{99}

\bibitem{ALMT14} Amelunxen, D., Lotz, M., McCoy, M.B., Tropp, J.A., Living on the edge: phase transitions in convex programs with random data. {\em Inf. Inference} {\bf 3} (2014), 224--294.

\bibitem{BV94} Baryshnikov, Y.M., Vitale, R.A., Regular simplices and Gaussian samples. {\em Discrete Comput. Geom.} {\bf 11} (1994), 141--147.

\bibitem{BCGTT19} Bonnet, G., Chasapis, G., Grote, J., Temesvari, D., Turchi, N., Threshold phenomena for high-dimensional random polytopes. {\em Commun. Contemp. Math.}  {\bf 21} (2019), no. 5, 1850038, 30 pp.

\bibitem{BKT19} Bonnet, G., Kabluchko, Z., Turchi, N., Phase transition for the volume of high-dimensional random polytopes. arXiv:1911.12696v1

\bibitem{BO19} Bonnet, G., O'Reilly, E., Facets of spherical random polytopes. arXiv:1908.04033 (2019).

\bibitem{CE67} Cover, T.M., Efron, B., Geometrical probability and random points on a hypersphere. {\em Ann. Math. Stat.} {\bf 38} (1967), 213--220.

\bibitem{DT09a} Donoho, D., Tanner, J., Observed universality of phase transitions in high-dimensional geometry, with implications for modern data analysis and signal processing. {\em Phil. Trans. R. Soc. A} {\bf 367} (2009), 4273--4293.

\bibitem{DT09b} Donoho, D., Tanner, J., Counting faces of randomly projected polytopes when the projection radically lowers dimension. {\em J. Amer. Math. Soc.} {\bf 22} (2009), 1--53.

\bibitem{DT10} Donoho, D., Tanner, J., Counting the faces of randomly-projected hypercubes and orthants, with applications. {\em Discrete Comput. Geom.} {\bf 43} (2010), 522--541.

\bibitem{DFD92} Dyer, M.E., F\"uredi, Z., McDiarmid, C., Volumes spanned by random points in the hypercube. {\em Random Structures Algorithms} {\bf 3} (1992), 91--106.


\bibitem{GG09} Gatzouras, D., Giannopoulos, A., Threshold for the volume spanned by random points with independent coordinates. {\em Israel J. Math.} {\bf 169} (2009), 125--153.

\bibitem{Gru68} Gr\"unbaum, B., Grassmann angles of convex polytopes. {\em Acta Math.} {\bf 121} (1968), 293--302.

\bibitem{HHRT15} H\"orrman, J., Hug, D., Reitzner, M., Th\"ale, C., Poisson polyhedra in high dimensions. {\em Adv. Math.} {\bf 281} (2015), 1--39.

\bibitem{HS16} Hug, D., Schneider, R., Random conical tessellations. {\em Discrete Comput. Geom.} {\bf 56} (2016), 395--426.

\bibitem{Klar} Klar, B., Bounds on tail probabilities of discrete distributions.
{\em  Probab.~Engrg.~Inform.~Sci.~} {\bf 14} (2000), 161--171.

\bibitem{Okamoto} Okamoto, M., Some inequalities relating to the partial sum of binomial probabilities.
{\em  Ann.~Inst.~Statist.~Math.} {\bf 10} (1958), 29--35.

\bibitem{Piv07} Pivovarov, P., Volume thresholds for Gaussian and spherical random polytopes and their duals. {\em Studia Math.} {\bf 183} (2007), 15--34.

\bibitem{SW08} Schneider, R., Weil, W., {\em Stochastic and Integral Geometry.} Springer, Berlin, 2008.

\bibitem{Shi06} Shiryaev, A.N., Probability--1. Third Edn., {\em Graduate Texts in Mathematics 95}, Springer, New York, 2006.

\bibitem{VS92} Vershik, A.M., Sporyshev, P.V., Asymptotic behavior of the number of faces of random polyhedra and the neighborliness problem. {\em Selecta Math. Soviet.} {\bf 11}, vol. 2 (1992), 181--201.

\end{thebibliography}
\end{document}